\documentclass[12pt]{amsart}
\usepackage{amssymb,amsmath}
\usepackage{graphpap}
\newtheorem{theorem}[subsection]{Theorem}
\newtheorem*{theorem*}{Theorem}

\newtheorem{proposition}[subsection]{Proposition}
\newtheorem{corollary}[subsection]{Corollary}

\theoremstyle{definition}
\newtheorem{definition}[subsection]{Definition}
\newtheorem{problem}[subsection]{Problem}
\newtheorem{example}[subsection]{Example}
\theoremstyle{remark}
\newtheorem{remark}[subsection]{Remark}
\newcommand{\mt}[1]{\operatorname{#1}}
\newcommand{\EEE}{{\mathbb E}}
\newcommand{\DDD}{{\mathbb D}}

\newcommand{\QQ}{{\mathbb Q}}

\newcommand{\CC}{{\mathbb C}}

\newcommand{\PP}{{\mathbb P}}

\newcommand{\Diff}{\mt{Diff}}
\newcommand{\Exc}{\mt{Exc}}
\newcommand{\codim}{\mt{codim}}

\title{On purely log terminal blow-ups }
\author{S.~A.~Kudryavtsev}
\date{}
\address{Department of Algebra, Faculty of Mathematics,
Moscow State Lomonosov University, 117234 Moscow,
Russia}
\email{kudryav@mech.math.msu.su}

\begin{document}
\begin{abstract} In this paper we prove the existence of purely log terminal
blow-up for Kawamata log terminal singularity and obtain the criterion for a
singularity to be weakly exceptional in terms of the exceptional divisor of
plt blow-up.
\end{abstract}
\maketitle

\section*{\bf {Introduction}}
The main aim of this note is to prove two results of the paper
\cite{Pr1} for non $\QQ$-factorial case.
The first one is the inductive blow-up existence theorem (theorem \ref{main1}) and
the second one is a criterion of weakly exceptionality (theorem \ref{main2}).
These blow-ups allow us to apply Shokurov's inductive method to the study of
singularities and in general case extremal contractions.
Using this method we can reduce the questions on structure, complementness and
exceptionality of singularity to a single exceptional divisor of purely
log terminal blow-up.
For any $\QQ$-factorial singularity a plt blow-up is the unique one that allows
to extend the complement of exceptional divisor to a global complement
(remark \ref{smain1}). For non $\QQ$-factorial klt singularity
such blow-ups differ from plt blow-ups by a small flopping contraction
(corollary \ref{smain2}).
In studying any $\QQ$-gorenstein singularities it is practically impossible to
select $\QQ$-factorial singularity class from the others.
That is why we have to apply the theorems and constructions which are true in the
general case. This paper also proves some results on the inductive method of
any lc singularity studies.

\subsubsection*{Acknowledgements}
I am grateful to Professor V.A.Iskovskikh and Professor Yu.G.Prokhorov for
useful discussions, criticisms and valuable remarks. The research was partially
supported by a grant 99-01-01132 from the Russian Foundation of Basic Research
and a grant INTAS-OPEN 97-2072.

\section{\bf Purely log terminal blow-ups and their properties}
All varieties are algebraic and are assumed to be defined over
$\CC$, complex number field. The results can be easily modified to the category
of analytic spaces. We use the terminology and notations of Log Minimal
Model Program and the main properties of complements given in
\cite{Koetal}, \cite{Pr1}, \cite{Pr2}. {\it A strictly lc singularity}
is called lc singularity, but is not klt singularity.

\begin{definition} \label{def1}
Let $X$ be a normal lc  variety and let $f\colon Y \to X$ be a blow-up such that
the exceptional locus of $f$ contains only one irreducible divisor
$E\ (\Exc(f)=E)$. Then $f:(Y,E) \to X$ is called  {\it a purely log terminal
(plt) blow-up}, if $K_Y+E$ is plt and
$-E$ is $f$-ample.
\end{definition}

\begin{remark} In the definition \ref{def1} it is demanded that divisor
$E$ must be $\QQ$-Cartier. Hence $Y$ is a $\QQ$-gorenstein variety.
\end{remark}

\begin{remark} \label{smain1}
\begin{enumerate}
\item If $X$ is klt then $-(K_Y+E)$ is $f$-ample.
Indeed, we have $K_Y+E=f^*K_X+(a(E,0)+1)E$ and $a(E,0)+1>0$.
\item If $X$ is strictly lc then $a(E,0)=-1$.
\item \cite[2.2]{Pr1} If $X$ is $\QQ$-factorial then $Y$ is also
$\QQ$-factorial and $\rho(Y/X)=1$. Hence, in definition
\ref{def1} for $\QQ$-factorial singularity it
is not necessary to demand that divisor $-E$ is $f$-ample,
because amplness takes place always. Note that every exceptional locus
component has codimension 1 for any birational contraction to
$\QQ$-factorial variety.
\item By inversion of adjunction $K_E+\Diff_E(0)$ is klt. If $K_E+\Diff_E(0)$ is
$n$-complementary then $K_Y+E$ is $n$-complementary and
$K_X$ is too \cite[2.8]{Pr1}.
\item \cite[2.2]{Pr1} Let $f_i:(Y_i,E_i)\to (X \ni P)$ be two
plt blow-ups. If $E_1$ and $E_2$ define the same discrete valuation of
function field $\mathcal{K}(X)$, then $f_1$ and $f_2$ are isomorphic.
\end{enumerate}
\end{remark}

\begin{problem} Describe the class of all weak log Del Pezzo surfaces,
generically $\PP^1$ and elliptic fibrations which can be exceptional divisors
of some plt blow-ups of a terminal, canonical,
$\varepsilon$-lt or lc singularities.
\end{problem}

The existence of plt blow-up for klt singularity follows from the next theorem.
\begin{theorem} \label{main1}
Let $X$ be  a klt variety and let $D\ne 0$ be a boundary on $X$ such that
$(X,D)$ is lc, but is not plt.
Suppose LMMP is true or $\dim X\le 3$.
Then there exists an inductive blow-up $f:Y\to X$
such that:
\begin{enumerate}
\item The exceptional locus of $f$ contains only one irreducible divisor
$E\ (\Exc(f)=E)$;
\item $K_Y+E+D_Y=f^*(K_X+D)$ is lc;
\item $K_Y+E+(1-\varepsilon)D_Y$ is plt and anti-ample over $X$
for any $\varepsilon>0$;
\item If $X$ is $\QQ$-factorial then $Y$ is also
$\QQ$-factorial and $\rho(Y/X)=1$.
\end{enumerate}
\begin{proof}
Let us consider the proof of this theorem for
$\QQ$-factorial singularities \cite[2.9]{Pr1}.
Let $g:\widehat Z \to X$ be a minimal log terminal modification
of $(X,D)$ and $\widehat E=\sum \widehat E_i$ be a reducible exceptional divisor
\cite[17.10]{Koetal}, \cite[9.1]{Sh2}.
By definition of such modification $\widehat Z$ is $\QQ$-factorial and
$K_{\widehat Z}+\widehat E+ D_{\widehat Z}=g^*(K_X+D)$ is dlt.
Since $X$ has only klt singularities then
$K_{\widehat Z}+\widehat E=g^*K_X+\sum (a(\widehat E_i,0)+1)\widehat E_i$
cannot be $g$-nef by numerical properties of contractions \cite[1.1]{Sh2}.
Run $K_{\widehat Z}+\widehat E$-MMP over $X$.
Hence at the last step we get a divisorial extremal contraction $g': \widetilde Z \to X$
(see diagram (1)) and
$K_{\widetilde Z}+{\widetilde E}+D_{\widetilde Z}={g'}^*(K_X+D)$ is lc, where
${\widetilde E}$ is an irreducible divisor.
Since $K_{\widetilde Z}+{\widetilde E}$ is plt then
$(\widetilde Z,\widetilde E+(1-\varepsilon)D_{\widetilde Z})$
is plt for any $\varepsilon>0$.
\par
If $X$ is $\QQ$-factorial then
$\Exc(g')=\widetilde E$, $\rho(\widetilde Z/X)=1$ and
$-\widetilde E$ is $g'$-ample by remark \ref{smain1}.
Therefore $g'$ is an inductive blow-up.

\begin{equation}
\begin{array}{ccccccc}
(\widehat Z,\ \widehat E=\sum \widehat E_i)& \dashrightarrow &
(\widetilde Z,\ \widetilde E)& \stackrel{\varphi}{\dashrightarrow} &
(Y',\ E')& \stackrel{\psi}{\longrightarrow} & (Y,E) \\
\Big\downarrow\vcenter{
\rlap{\tiny $g$}} & & \Big\downarrow\vcenter{
\rlap{\tiny $g'$}}
& &\Big\downarrow\vcenter{
\rlap{\tiny $f'$}}& &
\Big\downarrow\vcenter{
\rlap{\tiny $f$}}\\
X & = & X & = & X  & =& X \\
\end{array}
\end{equation}
\par
Assume that $\Exc(g')=\widetilde E \cup \Delta$ where $\Delta \ne \emptyset$ and
$\codim_{\widetilde Z}\Delta \ge 2$.
Obviously
$K_{\widetilde Z}+D_{\widetilde Z}\equiv -\widetilde E$ over $X$. Thus
$K_{\widetilde Z}+D_{\widetilde Z}$ is not
$g'$-nef and it is not negative for curves lying on $\widetilde E$.
Apply $K_{\widetilde Z}+D_{\widetilde Z}$-MMP. At the last step we get a
divisorial contraction
$f':Y' \to X$ and $\Exc(f')=E'$ is a irreducible divisor.
Note also that the birational map
$\varphi$ is a composition of log flips and
$K_{Y'}+E'+(1-\varepsilon)D_{Y'}$ is plt for any
$\varepsilon>0$. If $-E'$ is $f'$-ample then $f'$ is an inductive blow-up.
\par
Let $-E'$ is not $f'$-ample.
Since $-E'$ is $f'$-nef and
$K_{Y'}+E'={f'}^*K_X+(a(E',0)+1)E'$ where $a(E',0)+1>0$ then $-(K_{Y'}+E')$ is
$f'$-nef.
By Base Point Free Theorem \cite[3.1.2]{KMM} applied to klt divisor
$K_{Y'}+(1-\delta)E'$ ($0<\delta \ll 1$) the linear system
$|-n(K_{Y'}+E')|$ is free over $X$ for $n\gg 0$.
It gives small birational contraction
$\psi:(Y',E')\to (Y,E)$. Let $C$ be an exceptional curve. Since
$(K_{Y'}+E')\cdot C=0$ then $E'\cdot C=0$ and $K_{Y'}\cdot C=0$.
Therefore morphism
$\psi$ contracts the curve $C$ if and only if
$E'\cdot C=0$. Clearly,
the given blow-up $f:(Y,E)\to X$ is a required one.
\end{proof}
\end{theorem}

\begin{definition}\label{defexc}
Let $(X\ni P)$ be a lc singularity. It is said to be
{\it weakly exceptional} if there exists only one plt blow-up
(up to isomorphism). A lc
pair $(X,D)$ is said to be {\it exceptional}, where $D$ is boundary,
if there exists at most one divisor $E$ with discrepancy $a(E,D)=-1$.
A lc singularity $(X\ni P)$ is said to be {\it exceptional} if
$(X,D)$ is exceptional for any boundary $D$ whenever $K_X+D$ is lc.
\end{definition}

The LMMP is also used in the next corollary from theorem \ref{main1}.
\begin{corollary} \label{cor1}
Let $f:(Y,E)\to (X\ni P)$ be a plt blow-up of klt singularity and let
$\dim f(E)\ge 1$. Then there exists another plt blow-up of $(X\ni P)$.
Therefore the singularity is not weakly exceptional.
\begin{proof} Take two hyperplane sections $H_1$ and $H_2$ passing through the
point $P$ and not containing
$f(E)$. Let $c>0$ is a log canonical threshold of pair $(X,H_1+H_2)$.
Then $K_X+c(H_1+H_2)$ is not plt. The set
$f(E)$ is different from $LCS(X,c(H_1+H_2))$. Apply theorem
\ref{main1} for $(X,c(H_1+H_2))$. This completes the proof.
\end{proof}
\end{corollary}

\begin{proposition} \label{slcv}
Let $f:(Y,E)\to (X\ni P)$ be a plt blow-up of strictly lc singularity
and let
$\dim f(E)\ge 1$.  Then $(X\ni P)$ is not
exceptional singularity.
\begin{proof} As in proof of corollary \ref{cor1} there exists divisor
$D$ such that $(X,D)$ is lc, but is not plt and set
$f(E)$ is different from a minimal center $LCS(X,D)$.
Thus $(X\ni P)$ is not exceptional by definition.
\end{proof}
\end{proposition}

The LMMP is used in order to prove the necessary condition in the following
theorem.
\begin{theorem} \label{Slc}
Let $(X\ni P)$ be a strictly lc singularity. Then
\begin{enumerate}
\item If there exists a plt blow-up then it is the unique
(up to isomorphism).
\item The singularity is exceptional if and only if there exists a plt blow-up
$f\colon (Y,E)\to (X\ni P)$ such that $f(E)=P$.
\end{enumerate}
\begin{proof} The first statement follows from the properties (2) and (5)
of the remark \ref{smain1}. Let's prove the second part of theorem.
\par
{\it Necessity.} Assume that the singularity $(X\ni P)$ is exceptional. We will
construct a plt blow-up (cf. proof of theorem \ref{main1}).
Let $g':\widetilde Z \to X$ be a minimal log terminal modification of $X$ and
$\widetilde E=\sum \widetilde E_i$ be a reducible exceptional divisor. By definition
of such modification
$\widetilde Z$ is $\QQ$-factorial and
$K_{\widetilde Z}+\widetilde E={g'}^*K_X$ is dlt.
Since the singularity is exceptional then $\widetilde E$ is irreducible divisor
and $K_{\widetilde Z}+\widetilde E$ is plt.
Let $\Exc(g')=\widetilde E \cup \Delta$,
where $\Delta \ne \emptyset$ and
$\codim_{\widetilde Z}\Delta \ge 2$. Apply
$K_{\widetilde Z}$-MMP over $X$. Hence at the last step we get a divisorial
extremal contraction
$f':Y' \to X$ and $\Exc(f')=E'$ is an irreducible divisor. Divisor
$K_{Y'}+E'$ is also plt.
If $-E'$ is $f'$-ample then $f'$ is a required plt blow-up by proposition
\ref{slcv}.
Let $-E'$ is not $f'$-ample.
Since $-E'$ is $f'$-nef and
$K_{Y'}\equiv -E'$ over $X$, then  $K_{Y'}$ is
$f'$-nef.
By Base Point Free Theorem \cite[3.1.2]{KMM} the linear system
$|nK_{Y'}|$ is free over $X$ for $n\gg 0$. It gives small birational morphism
$h:(Y',E')\to (Y,E)$. The given blow-up
$f:(Y,E)\to X$ is plt because $K_Y$ is $f$-ample. By proposition
\ref{slcv} $f(E)=P$.
\par
{\it Sufficiency.} Conversely assume that there exists a required blow-up.
Note that
$E$ is an unique exceptional divisor with discrepancy
$a(E,0)=-1$. Let $(X,D)$ is any lc pair. Then $D=0$ because $f(E)=P$.
\end{proof}
\end{theorem}

\begin{corollary} Let $(X\ni P)$ be a strictly lc exceptional singularity. Then
the minimal index of complement is equal to gorenstein index of $(X\ni P)$.
\end{corollary}

\begin{remark} A minimal index of complementary is bounded for three dimensional
lc singularities \cite[7.1]{Sh1}. A hypothesis is that this index is not
more then 66.
For strictly lc exceptional singularities it was proved in papers
\cite{I} and \cite{Fuj}.
For non-exceptional non-isolated strictly lc singularities the gorenstein index
is not bounded \cite[5.1]{Fuj}.
\end{remark}

\begin{corollary}\cite[2.4]{PrI} Exceptional singularity is weakly exceptional.
\begin{proof} The existence of plt blow-up follows from theorems
\ref{main1} and \ref{Slc}. By
\cite[2.4]{PrI} such blow-up is unique.
\end{proof}
\end{corollary}

We have the next corollary by proofs of theorems \ref{main1} and \ref{Slc}.
\begin{corollary} \label{smain2} Notation as in definition
\ref{def1}. Assume that we don't require $-E$ to be ample
over $X$. Then such blow-up differs from a plt blow-up  by a
small flopping contraction.
\end{corollary}

\section{\bf Criterion of weakly exceptionality}
To prove $(3)\Rightarrow (1)$ in the next theorem we use LMMP.
\begin{theorem} \label{main2}
Let $(X\ni P)$ be a klt blow-up and let
$f\colon (Y,E)\to X$ be a plt blow-up of $P$.
Then the following conditions are equivalent:
\begin{enumerate}
\item
$(X\ni P)$ is not weakly exceptional;
\item
There is an effective  $\QQ$-divisor $D\ge\Diff_E(0)$ such that
$-(K_E+D)$ is ample and $(E,D)$ is not klt;
\item
There is an effective  $\QQ$-divisor $D\ge\Diff_E(0)$ such that
$-(K_E+D)$ is ample and $(E,D)$ is not lc.
\end{enumerate}
\begin{proof} The statements  $(1)\Rightarrow (2)$ and
$(2)\Rightarrow (3)$ follow from \cite[4.3]{Pr1}.
Let's prove
$(3)\Rightarrow (1)$. By corollary \ref{cor1} we can suppose $f(E)=P$.
It was proved in  \cite[theorem 4.3]{Pr1} the existence of effective
$\QQ$-Cartier divisor
$B=\sum b_iB_i$ such that $K_Y+E+B$ is lc, but is not plt.
Also $K_Y+E+B$ is anti-ample over $X$.
We can take very ample divisor
$H$ containing the minimal center of $LCS(Y,E+B)$.
There is a small rational number
$\varepsilon>0$ such that
$-(K_Y+E+B+\varepsilon H)$ is $f$-ample.
Replacing $B$ by $c(B+\varepsilon H)$ we can assume without loss of
generality that
$b_i<1$ for all $i$ ($c<1$ because $H$ contains a minimal center).
Denote $L=f(B)$.
Since $-(K_Y+E+B)$ is
$f$-ample and all $b_i<1$ then lc threshold $c'$ of pair $(X,L)$ is greater
than 1. If pair $(X,c'L)$ is plt then there is an effective
$\QQ$-Cartier divisor $L'$ that $(X,c'L+L')$ is lc, but is not plt.
By theorem \ref{main1} we have an inductive blow-up
$f':(Y',E')\to X$ of $(X,c'L+L')$. Moreover $K_Y+E+c'L_Y+L'_Y$ is not lc ($c'>1$).
Thus $f$ and $f'$ are not isomorphic plt blow-ups.
\end{proof}
\end{theorem}

\begin{example}\cite[4.7]{Pr1}, \cite[6.4]{Pr2} Two dimensional klt singularity
is weakly exceptional if and only if it has type
$\DDD_n, \EEE_6, \EEE_7$ or $\EEE_8$.
Among them the singularities of type $\DDD_n$ are not exceptional.
Two dimensional strictly lc singularity is weakly exceptional
(it is exceptional by theorem \ref{Slc}) if and only if it is simple elliptic or it has
type $\widetilde{\DDD}_4$ (see the minimal resolution graph in fig. 1),
$\widetilde{\EEE}_6, \widetilde{\EEE}_7, \widetilde{\EEE}_8$
(see the minimal resolution graph in fig. 2), where
$(n_1,n_2,n_3)=(3,3,3),(2,4,4)\  \text{or}\  (2,3,6)$ respectively.\\
\begin{picture}(140,70)(-60,0)
\put(30,35){\circle*{6}}
\put(10,35){\circle*{6}}
\put(50,35){\circle*{6}}
\put(30,55){\circle*{6}}
\put(30,15){\circle*{6}}

\put(30,35){\line(1,0){20}}
\put(10,35){\line(1,0){20}}
\put(30,35){\line(0,1){20}}
\put(30,15){\line(0,1){20}}

\put(32,37){\tiny{-n}}
\put(7,39){\tiny{-2}}
\put(47,39){\tiny{-2}}
\put(34,13){\tiny{-2}}
\put(27,59){\tiny{-2}}

\put(20,0){\tiny{Fig. 1}}

\put(150,50){\oval(40,16)}
\put(137,48){\tiny{$\CC^2/\mathbb Z_{n_1}$}}
\put(195,50){\circle*{6}}

\put(240,50){\oval(40,16)}
\put(227,48){\tiny{$\CC^2/\mathbb Z_{n_2}$}}

\put(170,50){\line(1,0){24}}
\put(195,50){\line(1,0){24}}

\put(195,22){\oval(40,16)}
\put(182,20){\tiny{$\CC^2/\mathbb Z_{n_3}$}}
\put(195,30){\line(0,1){20}}
\put(192,54){\tiny{-n}}

\put(185,0){\tiny{Fig. 2}}
\end{picture}
\end{example}

\begin{remark} A three dimensional terminal singularity is not
weakly exceptional \cite[4.8]{Pr1}.
\end{remark}

\begin{example} Let $(X\ni P)$ be an $n$-dimensional $(n\ge 3)$ canonical
hypersurface singularity given by equation
$(x^n_1+\cdots+x^n_n+x^{n+1}_{n+1}=0) \subset (\CC^{n+1},0)$.
The weighted blow-up of $\CC^{n+1}$ with weights $(n+1,\ldots,n+1,n)$
induces a plt blow-up of $P$. The obtained log Fano variety
$(E,\Diff_E(0))$ is
\begin{gather*}
(x^n_1+\cdots+x^n_n+x_{n+1}\subset
\PP(1,\ldots,1,n),\frac{n}{n+1}\{x_{n+1}=0\})=\\
=(\PP^{n-1},\frac{n}{n+1}Q_n),
\end{gather*}
where $Q_n$ is smooth hypersurface of degree $n$ in $\PP^{n-1}$. By theorem
\ref{main2} the singularity $(X\ni P)$ is weakly exceptional. The divisor
$\{x_{n+1}=0\}$ is 1-complement being not plt.
Therefore the singularity is not exceptional.
\end{example}


\begin{thebibliography}{99}

\bibitem{Fuj}
\emph{Fujino O.} The indices of log canonical
singularities,
e-print math.AG/9909035

\bibitem{I}
\emph{Ishii S.} The quotient of log--canonical
singularities by finite groups, appear in
Adv. Stud. in Pure Math., preprint TITECH-MATH 03-99 (\# 89)

\bibitem{KMM}
\emph{Kawamata Y., Matsuda K., Matsuki K.} Introduction to the minimal
model program, Adv. Stud. in Pure Math. 10 (1987), 283-360.

\bibitem{Koetal}
\emph{Kollar J. et al} Flips and abundance for algebraic
threefolds, Ast\'erisque 211,(1992).

\bibitem{Pr1}
\emph{Prokhorov Yu.~G.} Blow-ups of canonical singularities,
Algebra. Proc. Internat. Conf. on the Occasion of the
90th birthday of A.~G.~Kurosh, Moscow, Russia, May
25-30, 1998, Yu.~Bahturin ed., Walter de Gruyter, Berlin (2000), 301-317

\bibitem{Pr2}
\emph{Prokhorov Yu.~G.} Lectures on complements on log surfaces,
MSJ Memoirs 10, (2001), e-print math.AG/9912111

\bibitem{PrI}
\emph{Prokhorov Yu.~G., Ishii S.} Hypersurface exceptional
singularities, Intern. Journal of Math. 12 \textbf{6} (2001) 661--687
e-print math.AG/9910123

\bibitem{Sh2}
\emph{Shokurov V.~V.}
$3$-fold log flips, Russian Acad. Sci.
Izv. Math. \textbf{40} (1993) 93--202 \& \textbf{43} (1994)
527--558

\bibitem{Sh1}
\emph{Shokurov V.~V.} Complements on surfaces, J. of Math. Sci.
102 \textbf{2} (2000), 3876-3932. e-print math.AG/9711024

\end{thebibliography}
\end{document}